\documentclass[12pt]{amsart}
\usepackage{amssymb}
\usepackage{graphicx}
\usepackage{xcolor} 
\usepackage{tensor}

\usepackage[margin=1in, marginpar=.6in]{geometry} 

\usepackage{enumitem}

\allowdisplaybreaks				

\usepackage{mathtools}
\mathtoolsset{showonlyrefs}

\definecolor{green}{rgb}{0,0.8,0} 

\newtheorem{theorem}{Theorem}[section]

\newtheorem{lemma}[theorem]{Lemma}

\theoremstyle{definition}

\theoremstyle{remark}

\numberwithin{equation}{section}
\newcommand{\relphantom}[1]{\mathrel{\phantom{#1}}}

\makeatletter
\newcommand{\nrm}{\@ifstar{\nrmb}{\nrmi}}
\newcommand{\nrmi}[1]{\Vert{#1}\Vert}
\newcommand{\nrmb}[1]{\left\Vert{#1}\right\Vert}
\newcommand{\abs}{\@ifstar{\absb}{\absi}}
\newcommand{\absi}[1]{\vert{#1}\vert}
\newcommand{\absb}[1]{\left\vert{#1}\right\vert}
\newcommand{\brk}{\@ifstar{\brkb}{\brki}}
\newcommand{\brki}[1]{\langle{#1}\rangle}
\newcommand{\brkb}[1]{\left\langle{#1}\right\rangle}
\newcommand{\set}{\@ifstar{\setb}{\seti}}
\newcommand{\seti}[1]{\{#1\}}
\newcommand{\setb}[1]{\left\{ #1\right\}}
\makeatother

\newcommand{\VERT}[1]{{\left\vert\kern-0.25ex\left\vert\kern-0.25ex\left\vert #1 
    \right\vert\kern-0.25ex\right\vert\kern-0.25ex\right\vert}}

\DeclareMathOperator{\supp}{supp}

\newcommand{\aleq}{\lesssim}

\newcommand{\lap}{\Delta}

\newcommand{\ud}{\mathrm{d}}
\newcommand{\rd}{\partial}
\newcommand{\nb}{\nabla}

\newcommand{\bb}{\Big}

\newcommand{\lmb}{\lambda}

\newcommand{\bfu}{{\bf u}}
\newcommand{\bfv}{{\bf v}}

\newcommand{\bfB}{{\bf B}}

\newcommand{\bfJ}{{\bf J}}

\newcommand{\bfP}{{\bf P}}

\newcommand{\bfZ}{{\bf Z}}

\newcommand{\bfomg}{\boldsymbol{\omega}}

\newcommand{\bbR}{\mathbb R}

\newcommand{\bbT}{\mathbb T}

\newcommand{\calF}{\mathcal F}

\begin{document}

\title[]{Beale--Kato--Majda-type continuation criteria \\ for Hall- and electron-magnetohydrodynamics}
\author{Mimi Dai}%
\address{Department of Mathematics, Statistics and Computer Science, University of Illinois at Chicago, Chicago, IL 60607, USA}
\email{mdai@uic.edu} 

\author{Sung-Jin Oh}%
\address{Department of Mathematics, UC Berkeley, Berkeley, CA, USA 94720 and School of Mathematics, Korea Institute for Advanced Study, Seoul, Korea 02455}%
\email{sjoh@math.berkeley.edu}%


\begin{abstract}
We show that regular solutions to electron-MHD with resistivity can be continued as long as the time integral of the supremum of the current gradient remains finite. This dimensionless continuation criterion is analogous to the celebrated result of Beale--Kato--Majda for the incompressible Euler and Navier--Stokes equations. A similar continuation criterion, formulated in terms of the time integral of the supremum of the vorticity, velocity gradient and current gradient, is established for the Hall-MHD with resistivity as well.
\end{abstract}
\maketitle
\section{Introduction}
On $[0, \infty)_{t} \times \bbR^{3}_{x}$, the \emph{electron MHD equation (E-MHD)} with resistivity takes the form
\begin{equation} \label{eq:e-mhd}
\begin{aligned}
	\rd_{t} \bfB + \nb \times ( (\nb \times \bfB) \times \bfB) &= \lap \bfB, \\
	\nb \cdot \bfB &= 0
\end{aligned}
\end{equation}
where $\bfB(t) : \bbR^{3}_{x} \to \bbR^{3}$ for each $t \geq 0$.
The \emph{Hall MHD equation (H-MHD)} with resistivity takes the form 
\begin{equation} \label{eq:hall-mhd}
\begin{aligned}
	\rd_{t} \bfu + \bfu \cdot \nb \bfu - \bfB \cdot \nb \bfB + \nb p &= \nu \lap \bfu, \\
	\rd_{t} \bfB + \bfu \cdot \nb \bfB - \bfB \cdot \nb \bfu+ \nb \times ( (\nb \times \bfB) \times \bfB) &= \lap \bfB, \\
	\nb \cdot \bfu &=0,\\ 
	\nb \cdot \bfB &= 0
\end{aligned}
\end{equation}
where $\nu \geq 0$ is the constant viscosity parameter and unknowns $\bfu(t), \bfB(t) : \bbR^{3}_{x} \to \bbR^{3}$, $p(t): \bbR^{3}_{x} \to \bbR$ for each $t \geq 0$. Both systems \eqref{eq:e-mhd} and \eqref{eq:hall-mhd} are models from plasma physics, describing respectively the evolution of magnetic field $\bfB$ without or with the influence of the velocity field $\bfu$ of the ion flow. System \eqref{eq:e-mhd} is regarded as an approximating model of \eqref{eq:hall-mhd} when $\bfu$ is negligible and the constraint $- \bfB \cdot \nb \bfB + \nb p=0$ is relaxed. 
By vector calculus identities, the first equation of \eqref{eq:e-mhd} can be rewritten as
\begin{align*}
(\rd_{t}  - \lap) \bfB + \bfB \cdot \nb (\nb \times \bfB) - (\nb \times \bfB) \cdot \nb \bfB = 0,
\end{align*}
while the equation for $\rd_{t} \bfB$ in \eqref{eq:hall-mhd} can be rewritten as
\begin{align*}
(\rd_{t}  - \lap) \bfB + \bfu \cdot \nb \bfB - \bfB \cdot \nb \bfu + \bfB \cdot \nb (\nb \times \bfB) - (\nb \times \bfB) \cdot \nb \bfB = 0.
\end{align*}

The aim of this short note is to establish BKM-type continuation criteria for the above equations. The main results are as follows.
\begin{theorem}[BKM-type continuation criterion for E-MHD] \label{thm:cont-crit-J}
Let $\bfB$ be an $H^{\infty}$ solution to \eqref{eq:e-mhd} on $[0, T) \times \bbR^{3}_x$. If
\begin{equation} \label{eq:cont-crit-J}
	\int_{0}^{T} \nrm{\nb \bfJ(t)}_{L^{\infty}} \, \ud t < + \infty,
\end{equation}
where $\bfJ = \nb \times \bfB$, then $\bfB$ can be continued past $t = T$.
\end{theorem}

\begin{theorem}[BKM-type continuation criterion for Hall-MHD] \label{thm:cont-crit-H}
Let $(\bfu,\bfB)$ be an $H^{\infty}$ solution to \eqref{eq:hall-mhd} on $[0, T) \times \bbR^{3}_x$. If
\begin{equation} \label{eq:cont-crit-H}
	\int_{0}^{T} \nrm{\bfomg (t)}_{L^{\infty}}  + \nrm{\nb (\bfu - \bfJ)(t)}_{L^{\infty}} \, \ud t < + \infty,
\end{equation}
where $\bfomg = \nb \times \bfu$ and $\bfJ = \nb \times \bfB$, then $(\bfu,\bfB)$ can be continued past $t = T$.
\end{theorem}

Of course, the above theorems hold for sufficiently regular solutions, e.g., $H^{4}$ solutions instead of $H^{\infty}$. Moreover, the proofs presented below readily extend to solutions on $[0, T) \times \bbT^{3}$ instead of $[0, T) \times \bbR^{3}$. Some further remarks concerning these theorems are in order.

\smallskip \noindent {\it 1.} The main feature of these criteria is that they are formulated in terms of an $L^{1}_{t} L^{\infty}_{x}$-type norm involving the current $\bfJ = \nb \times \bfB$, similar to the celebrated criterion of Beale--Kato--Majda for incompressible Euler and Navier--Stokes equations \cite{BKM}. In the case of E-MHD, the quantity in \eqref{eq:cont-crit-J} is \emph{dimensionless}, in the sense that it is invariant under any scaling transformation\footnote{Note that the invariant scaling for the \eqref{eq:e-mhd} is the special case $\mu = \lmb^{2}$, but if we drop the resistivity $\lap u$ from the equation, then E-MHD is invariant under the two-parameter family of scaling transformations . We emphasize, however, that Theorem~\ref{thm:cont-crit-J} does \emph{not} hold in the absence of the term $\lap \bfB$; see 4.~below.} of the form $\bfB \mapsto \frac{\lmb^{2}}{\mu} \bfB(\mu^{-1} t, \lmb^{-1} x)$. Hall-MHD does not enjoy scaling invariance, but the quantity in \eqref{eq:cont-crit-H} is made out of $\nb \bfJ$ and quantities for $\bfu$ that scale like the original Beale--Kato--Majda criterion (i.e., $\bfomg$ and $\nb \bfu$ ).

\smallskip \noindent {\it 2.} Mathematically, a novel aspect of these criteria is that \eqref{eq:cont-crit-J} and \eqref{eq:cont-crit-H} impose very little control on the low frequency part of $\bfB$ due to the presence of several derivatives (i.e., $\nb \bfJ = \nb \nb \times \bfB$); see also the brief review of the previous results below. Indeed, obtaining a sufficiently strong control of $\bfB$ itself starting from \eqref{eq:cont-crit-J} and \eqref{eq:cont-crit-H} is the main difficulty of the proofs. We address this by making a (seemingly new) observation that \eqref{eq:cont-crit-J} or \eqref{eq:cont-crit-H} leads to a parabolic \emph{maximum principle} for $\abs{\bfB}^{2}$, and hence controls the growth of $\nrm{\bfB(t)}_{L^{\infty}}$. See Step~1 in Sections~\ref{sec:cont-crit-J} and \ref{sec:cont-crit-H}.
 
 \smallskip \noindent {\it 3.} Since $\bfJ$ is divergence-free, it is reasonable to suspect whether $\nb \bfJ$ in Theorem~\ref{thm:cont-crit-J} can be replaced by $\nb \times \bfJ$. Unfortunately, our proof breaks down, specifically the maximum principle for $\abs{\bfB}^{2}$. Same comment applies to the question of generalizing Theorem~\ref{thm:cont-crit-H} by replacing $\nb(\bfu - \bfJ)$ with $\bfomg - \nb \times \bfJ$.
 
\smallskip \noindent {\it 4.} As remarked earlier, \eqref{eq:cont-crit-J} is also invariant under the scalings of the non-resistive E-MHD, i.e., \eqref{eq:e-mhd} without $\lap \bfB$ on the right-hand side. One may therefore ask whether Theorem~\ref{thm:cont-crit-J} holds even in the non-resistive case. The answer is likely \emph{negative} in view of the recent illposedness results of Jeong and the second author \cite{JeOh3}, in the following sense. As is usual, our proof of Theorem~\ref{thm:cont-crit-J} proceeds by establishing \emph{persistence of regularity}, i.e., showing that a sufficiently high Sobolev norm of $\bfB$ (concretely, the $H^{4}$-norm) enjoys a uniform bound for $0 \leq t < T$ in terms of the control norm in \eqref{eq:cont-crit-J}. Since the time integral of the $H^{4}$-norm bounds the control norm in \eqref{eq:cont-crit-J}, such a persistence of regularity result would lead to a short time a-priori $H^{4}$-norm bound. However, this is impossible in view of the illposedness results in \cite{JeOh3}. A similar remark applies to Hall--MHD.

This failure is due to the fact that, in the absence of resistivity, the nonlinearity of $\nb \times (\bfB \times (\nb \times \bfB))$ in \eqref{eq:e-mhd} and \eqref{eq:hall-mhd} makes the $\bfB$-equation a quasilinear \emph{dispersive} equation, which is close to (say) a quasilinear Schr\"odinger equation but rather different from transport equations of incompressible fluid mechanics. We refer to \cite{JeOh1, JeOh2, JeOh3} for results in the non-resistive case based on this viewpoint.

\medskip
We end the introduction with a brief discussion of previously known continuation criteria for \eqref{eq:e-mhd} and \eqref{eq:hall-mhd}. The local well-posedness of classical solution of \eqref{eq:hall-mhd} (and \eqref{eq:e-mhd}) is known thanks to the work \cite{CDL} of Chae, Degond, and Liu. Moreover, the authors established the following continuation criterion for the local smooth solution beyond the time $T$,
\begin{equation}\label{crit-CDL}
\int_0^T \|\bfomg(t)\|_{\dot B^{0,\infty}_{\infty}}+\frac{1+\|\bfu(t)\|_{L^\infty}^2+\|\bfB(t)\|_{L^\infty}^2+\|\nabla\bfB(t)\|_{L^\infty}^2}{1+\log\left(1+\|\bfu(t)\|_{L^\infty}^2+\|\bfB(t)\|_{L^\infty}^2+\|\nabla\bfB(t)\|_{L^\infty}^2\right)} \, \ud t<\infty.
\end{equation}
A similar type of criterion as \eqref{crit-CDL} in term of Besov norms was obtained by Fan, Li and Nakamura \cite{FLN}. 
The Ladyzhenskaya-Prodi-Serrin type of criterion for \eqref{eq:hall-mhd} was established by Chae and Lee \cite{CL}, 
\begin{equation}\label{cl-criterion1}
\bfu\in L^q(0,T;L^p) \qquad \mbox { and } \qquad \nabla \bfB\in L^\gamma(0,T;L^\beta)
\end{equation}
with $\frac3p+\frac2q\leq 1$, $\frac3\beta+\frac2\gamma\leq 1$ and $p, \beta \in (3, \infty]$.
In the limit case $q=\gamma=2$, the authors also showed that if 
\begin{equation}\label{cl-bmo}
\int_0^T \|\bfu(t)\|_{BMO}^2+ \|\nabla\bfB(t)\|_{BMO}^2\, \ud t<\infty,
\end{equation}
the smooth solution can be continued beyond the time $T$. We note that the condition \eqref{cl-bmo} is weaker than \eqref{cl-criterion1} due to the embedding $L^\infty\hookrightarrow BMO$. The criterion \eqref{cl-bmo} was independently attained by Wan and Zhou \cite{WZ} where they obtained Osgood blowup criterion as well.

The first author established a regularity criterion for the 3D Hall-MHD \eqref{eq:hall-mhd} with $\nu=1$ on $\mathbb T^3$ in term of low frequencies of $\nabla\bfu$ and $\nabla\bfB$ in \cite{Dai-hmhd}. The author defined time dependent dissipation wavenumber $\kappa_{\bfu}(t)$ and $\kappa_{\bfB}(t)$ with respect to $\bfu$ and $\bfB$. 
Assuming 
\begin{equation}\label{crit-low}
\int_0^T\left(\|P_{<\kappa_{\bfu}(t)}\nabla \bfu(t)\|_{ B^{0,\infty}_{\infty}}+2^{\kappa_{\bfB}(t)}\|P_{<\kappa_{\bfB}(t)}\nabla \bfB(t)\|_{ B^{0,\infty}_{\infty}}\right)\, \ud t<\infty,
\end{equation}
the local smooth solution $(\bfu(t), \bfB(t))$ can be extended over $T$. The projection notation $P_{<}$ is introduced in Section \ref{sec:prelim}. 
It was shown in \cite{Dai-hmhd} that the low frequency condition \eqref{crit-low} is weaker than \eqref{cl-criterion1} and \eqref{cl-bmo}; nevertheless, the new criterion \eqref{eq:cont-crit-H} does \emph{not} imply \eqref{crit-low}, and vice versa. We also mention the articles \cite{HeAhHaZh, Ye}, in which a continuation criterion for \eqref{eq:hall-mhd} with\footnote{We note that these proofs rely heavily on the fact that the variable $\bfZ := \bfB + \bfomg$ solves an advection--diffusion equation (with the advection given by $\bfu$), which holds only when $\nu = 1$ (i.e., viscosity and resistivity match).} $\nu = 1$ only in terms of the $\bfu$ was derived.

Finally, considering the E-MHD without resistivity, i.e. \eqref{eq:e-mhd} with $\Delta\bfB$ removed, 
the first author, Guerra and Wu \cite{DGW} recently proved that a class of self-similar solutions do not blow up at time $T$ if the BKM-type condition 
\begin{equation}\label{crit-self-s}
\lim_{t\to T^-}\int_0^t\|\nabla\times \bfJ(\tau)\|_{L^\infty}\, \ud \tau<\infty
\end{equation} 
holds. Despite the above-mentioned issue with the BKM-type criterion \eqref{crit-self-s} for E-MHD without resistivity in general, the work above shows such continuation criterion for a class of self-similar solutions.

\subsection*{Acknowledgement}
M. Dai is grateful for the support from Simons Foundation and the National Science Foundation through the grants NSF-DMS-2009422 and NSF-DMS-2308208.
S.-J.~Oh would like to thank the hospitality of University of Chicago Illinois, where this work was carried out. S.-J.~Oh was partially supported by a Sloan Research Fellowship and a National Science Foundation CAREER Grant under NSF-DMS-1945615.

\bigskip

\section{Preliminaries} \label{sec:prelim}
We follow the standard convention of using $C$ to denote positive constants that may vary from expression to expression. We also use the shorthand $A \aleq B$ for the inequality $A \leq C B$ for some $C > 0$.

We employ the following convention for Littlewood--Paley projections. Fix the function $m_{<0} \in C^{\infty}_{c}(\bbR^{3})$ that is radial and non-increasing (in $r$) with $\supp m_{<0} \subseteq \set{\xi \in \bbR^{3} : \abs{\xi} < 2}$ and $m_{<0}(\xi) = 1$ on $\set{\xi \in \bbR^{3} : \abs{\xi} < 1}$. We define $P_{< k} f := \calF^{-1}[m_{<0}(2^{-k} \xi) \calF[f]]$, and $P_{k} f := P_{< k+1} f - P_{<k} f$  (where $\calF$ is the Fourier transform). 

We will be using the (inhomogeneous) Besov norm $\nrm{\cdot}_{B^{0, \infty}_{\infty}}$ (with regularity, integrability and summability exponents $0, \infty, \infty$), which is defined as
\begin{equation*}
\nrm{f}_{B^{0, \infty}_{\infty}} = \sup_{k \in \set{0, 1, 2, \ldots}} \nrm{P_{k} f}_{L^{\infty}} + \nrm{P_{<0} f}_{L^{\infty}}.
\end{equation*}
Some basic facts about $B^{0, \infty}_{\infty}$ are as follows. By the $L^{\infty}$-boundedness of $P_{k}$, we have (for every $f \in L^{\infty}$)
\begin{equation} \label{eq:Linfty-besov}
	\nrm{f}_{B^{0, \infty}_{\infty}} \aleq \nrm{f}_{L^{\infty}}.
\end{equation}
Since $P_{k} (-\lap)^{-1} \nb \times (\nb \times \cdot)$ is bounded for each $k \geq 0$ (uniformly in $k$), we have (as long as both sides are finite)
\begin{equation} \label{eq:besov-sio}
	\nrm{\nb \bfv}_{B^{0, \infty}_{\infty}} \aleq \nrm{\nb \times \bfv}_{B^{0, \infty}_{\infty}}.
\end{equation}

The following simple interpolation bound will be useful.
\begin{lemma} \label{lem:Linfty-interpolate}
Let $\bfv \in H^{4}(\bbR^{3})$. We have
\begin{equation} \label{lem:Linfty-interpolate}
	\nrm{\nb \bfv}_{L^{\infty}} \aleq \nrm{\bfv}_{L^{\infty}}^{\frac{1}{2}} \nrm{\nb (\nb \times \bfv)}_{B^{0, \infty}_{\infty}}^{\frac{1}{2}}.
\end{equation}
\end{lemma}
The regularity condition $\bfv \in H^{4}$ is inessential and has only been fixed for concreteness. Indeed, by an approximation argument, this lemma may be extended to $\bfv$ belonging to a larger space (which we will not detail).
\begin{proof}
Let $k_{0}$ be a nonnegative integer to be fixed later. We split $\nb \bfv = \sum_{k = k_{0}}^{\infty} P_{k} \nb \bfv + P_{<k_{0}} \nb \bfv$ and bound
\begin{align*}
\nrm{\nb \bfv}_{L^{\infty}} &\aleq \sum_{k=k_{0}} \nrm{P_{k} \nb \bfv}_{L^{\infty}} + \nrm{P_{<k_{0}} \nb \bfv}_{L^{\infty}} \\
&\aleq \sum_{k=k_{0}} 2^{-k} \nrm{P_{k} \nb (\nb \times \bfv)}_{B^{0, \infty}_{\infty}} + 2^{k_{0}} \nrm{P_{<k_{0}}\bfv}_{L^{\infty}} \\
&\aleq 2^{-k_{0}} \nrm{\nb (\nb \times \bfv)}_{B^{0, \infty}_{\infty}} + 2^{k_{0}} \nrm{\bfv}_{L^{\infty}}.
\end{align*}
where we used, among others, Bernstein's inequality and $L^{\infty}$-boundedness of $P_{k} (-\lap)^{-1} \nb \times (\nb \times \cdot)$. Now optimizing the choice of $k_{0}$ such that $2^{k_0}\sim \|\bfv\|_{L^\infty}^{-\frac12} \nrm{\nb (\nb \times \bfv)}_{B^{0, \infty}_{\infty}}^{\frac12}$, we obtain the lemma. \qedhere
\end{proof}

We will also need the following BKM-type singular integral bounds.
\begin{lemma} [BKM-type bound] \label{lem:bkm-sio}
For $\bfv \in H^{4} \cap L^{p}$ ($1 \leq p \leq \infty$), we have
\begin{align} 
	\nrm{\nb \bfv(t)}_{L^{\infty}}
	& \aleq \nrm{\nb \times \bfv(t)}_{B^{0, \infty}_{\infty}} \log(2 + \nrm{\nb^{4} \bfv(t)}_{L^{2}}) + \nrm{\bfv}_{L^{p}}, \label{eq:bkm-sio-1}	\\
	\nrm{\nb^{2} \bfv(t)}_{L^{\infty}}
	& \aleq \nrm{\nb (\nb \times \bfv)(t)}_{B^{0, \infty}_{\infty}} \log(2 + \nrm{\nb^{4} \bfv(t)}_{L^{2}}) + \nrm{\bfv}_{L^{p}}.\label{eq:bkm-sio-2}
\end{align}
\end{lemma}
As before, the regularity condition $\bfv \in H^{4}$ has only been fixed for concreteness; with essentially the same proof as below, $H^{4}$ may be replaced by any $H^{s}$ with $s > \frac{3}{2}+1$ and $s > \frac{3}{2} + 2$ in \eqref{eq:bkm-sio-1} and \eqref{eq:bkm-sio-2}, respectively.
\begin{proof}
Fix $m \in \set{1, 2}$ and let $k_{0}$ be a nonnegative integer to be fixed later. We begin with the decomposition
\begin{equation*}
	\nb^{m} \bfv = \sum_{k = k_{0}}^{\infty} P_{k} \nb^{m} \bfv + \sum_{k=0}^{k_{0}-1} P_{k} \nb^{m} \bfv + P_{< 0} \nb^{m} \bfv.
\end{equation*}
We may estimate
\begin{align*}
	\nrm{\nb^{m} \bfv}_{L^{\infty}} 
	&\leq \sum_{k = k_{0}}^{\infty} \nrm{P_{k} \nb^{m} \bfv}_{L^{\infty}} + \sum_{k = 0}^{k_{0}-1} \nrm{P_{k} \nb^{m} \bfv}_{L^{\infty}} + \nrm{P_{< 0} \nb^{m} \bfv}_{L^{\infty}} \\
	&\aleq \sum_{k = k_{0}}^{\infty} 2^{-(4 - \frac{3}{2} - m) k} \nrm{P_{k} \nb^{4} \bfv}_{L^{2}}
	+ \sum_{k=0}^{k_{0}-1} \nrm{P_{k} \nb^{m-1} (\nb \times \bfv)}_{L^{\infty}}
	+ \nrm{\bfv}_{L^{p}} \\
	&\aleq 2^{-(4 - \frac{3}{2} - m) k_{0}} \nrm{\nb^{4} \bfv}_{L^{2}} + k_{0} \nrm{\nb^{m-1} (\nb \times \bfv)}_{B^{0, \infty}_{\infty}} + \nrm{\bfv}_{L^{p}},
\end{align*}
where we used, among others, Bernstein's inequality and $L^{\infty}$-boundedness of $P_{k} (-\lap)^{-1} \nb \times (\nb \times \cdot)$. Note that $4 - \frac{3}{2} - m > 0$ for $m = 1, 2$. Now optimizing the choice of $k_{0}$, we obtain the lemma.
\end{proof}

For a comprehensive introduction to Littlewood--Paley theory and harmonic analysis background for this paper, we refer the reader to \cite{BCD}.

\bigskip

\section{Proof of Theorem~\ref{thm:cont-crit-J}} \label{sec:cont-crit-J}
In view of the known local well-posedness results (as well as the conservation of $\nrm{\bfB}_{L^{2}}$) \cite{CDL}, it suffices to show that 
\begin{equation*}
\limsup_{t \to T} \nrm{\nb^{4} \bfB(t)}_{L^{2}}  < + \infty.
\end{equation*}

\smallskip
\noindent{\it Step~1. Maximum principle for $\abs{\bfB}^{2}$.}
Observe that the scalar quantity $\abs{\bfB}^{2}$ obeys the following equation:
\begin{equation*}
	(\rd_{t} - \lap - \bfJ \cdot \nb) \abs{\bfB}^{2} 
	= - 2 \bfB \cdot \left( (\bfB \cdot \nb) \bfJ \right) - 2 \abs{\nb \bfB}^{2},
\end{equation*}
where we remind the reader that $\bfJ = \nb \times \bfB$.
Let $K > 0$ be a constant to be determined later, and define
\begin{equation*}
U(t, x) := \exp\left( - K I(t) \right) \abs{\bfB}^{2}(t, x), \quad I(t) := \int_{0}^{t} \nrm{\nb \bfJ(t')}_{L^{\infty}} \, \ud t'.
\end{equation*}
Note that $U$ obeys the following:
\begin{align*}
	 (\rd_{t} - \lap - (\nb \times \bfB) \cdot \nb) U
	&= \exp\left( - K I(t) \right) \bb[- K \nrm{\nb \bfJ(t)}_{L^{\infty}} \abs{\bfB}^{2}(t,x) \\
	&\relphantom{= \exp\left( - K I(t) \right) \bb[}
	+ (\rd_{t} - \lap - (\nb \times \bfB) \cdot \nb) \abs{\bfB}^{2}(t,x) \bb] \\
	&= \exp\left( - K I(t) \right)\bb[- K \nrm{\nb \bfJ(t)}_{L^{\infty}} \abs{\bfB}^{2}(t,x) \\
	&\relphantom{= \exp\left( - K I(t) \right) \bb[}
	- 2 \bfB \cdot \left( (\bfB \cdot \nb) \nb \times \bfB \right) 
	 - 2 \abs{\nb \bfB}^{2} \bb].
\end{align*}
Taking $K$ sufficiently large, we may arrange so that
\begin{align*}
 (\rd_{t} - \lap - (\nb \times \bfB) \cdot \nb) U \leq 0.
\end{align*}
At this point, we may apply the maximum principle to conclude that $U$ is uniformly bounded by $\nrm{U(0)}_{L^{\infty}}$ up to $t = T$. As a consequence,
\begin{equation} \label{eq:cont-crit-J:pf:0}
	\nrm{\bfB}_{L^{\infty}L^{\infty}([0, T) \times \bbR^{3})} \leq \exp \left(C \int_{0}^{T} \nrm{\nb \bfJ(t)}_{L^{\infty}} \, \ud t \right) \nrm{\bfB(0)}_{L^{\infty}}.
\end{equation}

\smallskip
\noindent{\it Step~2. Singular integral bounds.}
We claim that
\begin{align} 
	\nrm{\nb \bfB}_{L^{\infty}} &\aleq \nrm{\bfB}_{L^{\infty}}^{\frac{1}{2}} \nrm{\nb \bfJ}_{L^{\infty}}^{\frac{1}{2}}, \label{eq:cont-crit-J:pf:nbB} \\
	\nrm{\nb^{2} \bfB}_{L^{\infty}} &\aleq \nrm{\nb \bfJ}_{L^{\infty}} \log(2 + \nrm{\nb^{4} \bfB}_{L^{2}}) + \nrm{\bfB}_{L^{\infty}}. \label{eq:cont-crit-J:pf:nbnbB}
\end{align}
The first bound follows from Lemma~\ref{lem:Linfty-interpolate} and \eqref{eq:Linfty-besov}. The second bound follows from \eqref{eq:bkm-sio-2} (from Lemma~\ref{lem:bkm-sio}) with $p = \infty$ and \eqref{eq:Linfty-besov}.

\smallskip
\noindent{\it Step~3. Propagation of $\nrm{\nb^{4} \bfB}_{L^{2}}$.}
Let $I$ be a multi-index with $\abs{I} = 4$. Taking $\rd^{I}$ of the first equation in \eqref{eq:e-mhd}, we obtain
\begin{align*}
	(\rd_{t} - \lap) \rd^{I} \bfB  + \nb \times ((\nb \times \rd^{I} \bfB) \times \bfB)
	&= \sum_{I', I'' : I' + I'' = I, \, \abs{I'} \leq 3} c_{I', I''} \nb \times ((\nb \times \rd^{I'} \bfB) \times \rd^{I''} \bfB)
\end{align*}
where each $c_{I', I''}$ is some combinatorial coefficient. Taking the inner product of the equation above with $\rd^{I} \bfB$, integrating on $\bbR^{3}$ and summing over all $I$ with $\abs{I} = 4$, we obtain
\begin{align*}
	\frac{1}{2} \frac{\ud}{\ud t} \nrm{\nb^{4} \bfB}_{L^{2}}^{2} + \nrm{\nb^{5} \bfB}_{L^{2}}^{2} &\aleq \sum_{I', I'' : \abs{I'} + \abs{I''} = 4, \, \abs{I'} \leq 3} \int \abs{\nb \times ((\nb \times \rd^{I'} \bfB) \times \rd^{I''} \bfB)} \abs{\nb^{4} \bfB} \, \ud x \\
	&\aleq \sum_{m_{1} + m_{2} = 6, \, 1 \leq m_{1} \leq m_{2} \leq 5} \int \abs{\nb^{m_{1}} \bfB} \abs{\nb^{m_{2}} \bfB} \abs{\nb^{4} \bfB} \, \ud x \\
\end{align*}
We estimate each summand on the last line as follows.
\begin{itemize}
\item {\bf Case~1: $m_{1} = 1$ and $m_{2} = 5$.}
Let $\eta > 0$ be a parameter to be chosen later. By Cauchy--Schwarz and \eqref{eq:cont-crit-J:pf:nbB}, we have
\begin{align*}
 \int \abs{\nb \bfB} \abs{\nb^{5} \bfB} \abs{\nb^{4} \bfB} \, \ud x 
 &\leq \eta \nrm{\nb^{5} \bfB}_{L^{2}}^{2} + \eta^{-1} \nrm{\nb \bfB}_{L^{\infty}}^{2} \nrm{\nb^{4} \bfB}_{L^{2}}^{2} \\
& \leq \eta \nrm{\nb^{5} \bfB}_{L^{2}}^{2} + C \eta^{-1} \nrm{\bfB}_{L^{\infty}} \nrm{\nb \bfJ}_{L^{\infty}} \nrm{\nb^{4} \bfB}_{L^{2}}^{2}.
\end{align*}
\item {\bf Case~2: $2 \leq m_{1} \leq m_{2} \leq 4$.}
In this case, by interpolation and \eqref{eq:cont-crit-J:pf:nbnbB}, we may estimate
\begin{align*}
\int \abs{\nb^{m_{1}} \bfB} \abs{\nb^{m_{2}} \bfB} \abs{\nb^{4} \bfB} \, \ud x
&\aleq \nrm{\nb^{2} \bfB}_{L^{\infty}} \nrm{\nb^{4} \bfB}_{L^{2}}^{2} \\
&\aleq \left( \nrm{\nb \bfJ}_{L^{\infty}} \log (2 + \nrm{\nb^{4} \bfB}_{L^{2}}) + \nrm{\bfB}_{L^{\infty}} \right) \nrm{\nb^{4} \bfB}_{L^{2}}^{2}.
\end{align*}
\end{itemize}
In conclusion, we have
\begin{align*}
	&\frac{1}{2} \frac{\ud}{\ud t} \nrm{\nb^{4} \bfB}_{L^{2}}^{2} + \nrm{\nb^{5} \bfB}_{L^{2}}^{2} \\
	&\aleq \eta \nrm{\nb^{5} \bfB}_{L^{2}}^{2} + \left[\left(\log (2 + \nrm{\nb^{4} \bfB}_{L^{2}}) + (1+\eta^{-1}) \nrm{\bfB}_{L^{\infty}} \right) \nrm{\nb \bfJ}_{L^{\infty}} + \nrm{\bfB}_{L^{\infty}} \right]\nrm{\nb^{4} \bfB}_{L^{2}}^{2}.
\end{align*}
Choosing $\eta > 0$ sufficiently small, we may absorb the first term on the right-hand side into the left-hand side. Recall also that we have proved the boundedness of $\nrm{\bfB}_{L^{\infty}}$ in Step~1. Thus, by generalized Gr\"onwall's inequality (as in the usual proof of the Beale--Kato--Majda theorem), it follows that $\nrm{\nb^{4} \bfB}_{L^{2}}$ is uniformly bounded up to $t = T$. \hfill \qedsymbol

\bigskip

\section{Proof of Theorem~\ref{thm:cont-crit-H}} \label{sec:cont-crit-H}
As in the case of E-MHD, in view of the known local well-posedness results (as well as the conservation of $\nrm{\bfu}^{2} + \nrm{\bfB}_{L^{2}}^{2}$) such as \cite{CDL}, it suffices to show that 
\begin{equation*}
\limsup_{t \to T} \left( \nrm{\nb^{4} \bfu(t)}_{L^{2}} + \nrm{\nb^{4} \bfB(t)}_{L^{2}} \right)  < + \infty.
\end{equation*}

\noindent{\it Step~1. Maximum principle for $\abs{\bfB}^{2}$.}
We first show the maximum principle,
\begin{equation}\label{max-uJ} 
	\nrm{\bfB}_{L^{\infty} L^{\infty}([0, T) \times \bbR^{3})} \leq \exp \left(C \int_{0}^{T} \nrm{\nb (\bfu - \bfJ)(t)}_{L^{\infty}} \, \ud t \right) \nrm{\bfB(0)}_{L^{\infty}}.
\end{equation}
Indeed, taking the dot product of the second equation of \eqref{eq:hall-mhd} with $\bfB$ yields 
\begin{align*}
	(\rd_{t} - \lap + (\bfu - \bfJ) \cdot \nb) \abs{\bfB}^{2}
	&= 2 \bfB \cdot \left( (\bfB \cdot \nb) (\bfu - \bfJ) \right) - 2 \abs{\nb \bfB}^{2}.
\end{align*}
As in the E-MHD case, we define
\begin{equation*}
U(t, x) := \exp\left( - K I(t) \right) \abs{\bfB}^{2}(t, x), \quad I(t) := \int_{0}^{t}\nrm{\nb (\bfu - \bfJ)(t')}_{L^{\infty}}\, \ud t'
\end{equation*}
for large enough $K>0$.
A straightforward computation shows that $U$ satisfies
\begin{align*}
	 &(\rd_{t} - \lap + (\bfu- \bfJ) \cdot \nb) U\\
	=& \exp\left( - K I(t) \right) \bb[- K \nrm{\nb (\bfu - \bfJ)(t)}_{L^{\infty}} \abs{\bfB}^{2}(t,x) \\
	&\relphantom{= \exp\left( - K I(t) \right) \bb[}
	+ (\rd_{t} - \lap + (\bfu- \bfJ) \cdot \nb) \abs{\bfB}^{2}(t,x) \bb] \\
	=& \exp\left( - K I(t) \right)\bb[- K\nrm{\nb (\bfu - \bfJ)(t)}_{L^{\infty}} \abs{\bfB}^{2}(t,x) \\
	&\relphantom{= \exp\left( - K I(t) \right) \bb[}
	+ 2 \bfB \cdot \left( (\bfB \cdot \nb) (\bfu- \bfJ) \right) 
	 - 2 \abs{\nb \bfB}^{2} \bb].
\end{align*}
In view of the right hand side of the equation above, we can choose $K$ sufficiently large such that
\begin{align*}
 (\rd_{t} - \lap + (\bfu- \bfJ) \cdot \nb) U \leq 0.
\end{align*}
Hence applying the maximum principle to $U$ implies \eqref{max-uJ}.

\smallskip
\noindent{\it Step~2. Singular integral bounds.}
We claim that
\begin{align} 
	\nrm{\nb \bfu}_{L^{\infty}} &\aleq \nrm{\bfomg}_{L^{\infty}} \log(2 + \nrm{\nb^{4} \bfu}_{L^{2}}) + \nrm{\bfu}_{L^{2}}, \label{eq:cont-crit-H:pf:nbu} \\
	\nrm{\nb \bfB}_{L^{\infty}} &\aleq \nrm{\bfB}_{L^{\infty}}^{\frac{1}{2}} (\nrm{\bfomg}_{L^{\infty}}  + \nrm{\nb(\bfu - \bfJ)}_{L^{\infty}})^{\frac{1}{2}}, \label{eq:cont-crit-H:pf:nbB} \\
	\nrm{\nb^{2} \bfB}_{L^{\infty}} &\aleq (\nrm{\bfomg}_{L^{\infty}} + \nrm{\nb(\bfu - \bfJ)}_{L^{\infty}}) \log(2 + \nrm{\nb^{4} \bfB}_{L^{2}}) + \nrm{\bfB}_{L^{\infty}}. \label{eq:cont-crit-H:pf:nbnbB}
\end{align}
Indeed, \eqref{eq:cont-crit-H:pf:nbu} follows from \eqref{eq:bkm-sio-1} (in Lemma~\ref{lem:bkm-sio}) with $\bfv = \bfu$ and $p = 2$, and \eqref{eq:Linfty-besov}. To prove \eqref{eq:cont-crit-H:pf:nbB}, we use Lemma~\ref{lem:Linfty-interpolate}, \eqref{eq:Linfty-besov} and \eqref{eq:besov-sio} to bound
\begin{align*}
\nrm{\nb \bfB}_{L^{\infty}}
&\aleq \nrm{\bfB}_{L^{\infty}}^{\frac{1}{2}} \nrm{\nb \bfJ}_{B^{0, \infty}_{\infty}}^{\frac{1}{2}} \\
& \aleq \nrm{\bfB}_{L^{\infty}}^{\frac{1}{2}} (\nrm{\nb (\bfu - \bfJ)}_{B^{0, \infty}_{\infty}} + \nrm{\nb \bfu}_{B^{0, \infty}_{\infty}})^{\frac{1}{2}} \\
& \aleq \nrm{\bfB}_{L^{\infty}}^{\frac{1}{2}}(\nrm{\nb (\bfu - \bfJ)}_{B^{0, \infty}_{\infty}} + \nrm{\bfomg}_{B^{0, \infty}_{\infty}})^{\frac{1}{2}}.
\end{align*}
Finally, \eqref{eq:cont-crit-H:pf:nbnbB} follows from \eqref{eq:bkm-sio-1} (in Lemma~\ref{lem:bkm-sio}) with $\bfv = \bfu$ and $p = 2$, and a similar argument as above to estimate $\nrm{\nb \bfJ}_{B^{0, \infty}_{\infty}}$.

\smallskip
\noindent{\it Step~3. Propagation of $\nrm{\nb^{4} \bfu}_{L^{2}}$ and $\nrm{\nb^{4} \bfB}_{L^{2}}$.}
Again for the multi-index $I$ with $|I|=4$ we have
\begin{equation}\label{eq-h4}
\begin{split}
(\partial_t-\nu\Delta)\partial^{I}\bfu=&-\nabla\partial^I\bfP-\partial^I(\bfu\cdot\nabla\bfu)+\partial^I(\bfB\cdot\nabla\bfB),\\
(\partial_t-\Delta)\partial^{I}\bfB=&-\partial^I(\bfu\cdot\nabla\bfB)+\partial^I(\bfB\cdot\nabla\bfu)-\partial^I\nabla\times((\nabla\times\bfB)\times\bfB).
\end{split}
\end{equation}
Taking the inner product of the first equation of \eqref{eq-h4} with $\partial^I\bfu$ and the second equation with $\partial^I\bfB$, and adding the two resulted equations gives
\begin{equation}\label{est-h4-uB-2}
\begin{split}
&\frac12\frac{\ud}{\ud t}\left(\|\nabla^4\bfu\|_{L^2}^2+ \|\nabla^4\bfB\|_{L^2}^2\right)+\nu \|\nabla^5\bfu\|_{L^2}^2+\|\nabla^5\bfB\|_{L^2}^2\\
=&-\int \partial^I(\bfu\cdot\nabla\bfu)\cdot\partial^I\bfu\,\ud x+\int \partial^I(\bfB\cdot\nabla\bfB)\cdot\partial^I\bfu\,\ud x-\int \partial^I(\bfu\cdot\nabla\bfB)\cdot\partial^I\bfB\,\ud x\\
&+\int \partial^I(\bfB\cdot\nabla\bfu)\cdot\partial^I\bfB\,\ud x-\int \partial^I\nabla\times((\nabla\times\bfB)\times\bfB)\cdot\partial^I\bfB\,\ud x.
\end{split}
\end{equation}
For the first integral on the right hand side of \eqref{est-h4-uB-2}, we have
\begin{equation}\notag
\int \partial^I(\bfu\cdot\nabla\bfu)\cdot\partial^I\bfu\,\ud x=\sum_{I'+I''=3, |I''|\leq 3}c_{I',I''} \int (\partial^{I'}\bfu\cdot\nabla)\partial^{I''}\bfu\cdot\partial^I\bfu \,\ud x
\end{equation}
where we employed the cancellation
\[ \int (\bfu\cdot\nabla) \partial^I\bfu\cdot\partial^I\bfu\,\ud x =0.\]
Therefore,
\begin{equation}\notag
\left|\int \partial^I(\bfu\cdot\nabla\bfu)\cdot\partial^I\bfu\,\ud x \right|\aleq\sum_{m_1+m_2=5, 1\leq m_1\leq m_2\leq 4}\int |\nabla^{m_1} \bfu||\nabla^{m_2} \bfu||\nabla^4 \bfu|\, \ud x.
\end{equation}
In the case of $m_1=1$ and $m_2=4$, we obtain
\begin{equation}\notag
\begin{split}
\int |\nabla\bfu||\nabla^4 \bfu||\nabla^4 \bfu|\, \ud x\aleq& \|\nabla\bfu\|_{L^\infty}\|\nabla^4\bfu\|_{L^2}^2\\
\aleq& \left(1+\|\bfomg\|_{L^\infty}\log(2+\|\nabla^4\bfu\|_{L^2}) + \nrm{\bfu}_{L^{2}}\right)  \|\nabla^4\bfu\|_{L^2}^2
\end{split}
\end{equation}
using \eqref{eq:cont-crit-H:pf:nbu}. In the case of $2\leq m_1\leq m_2\leq 3$, we apply H\"older's inequality and Gagliardo-Nirenberg interpolation inequality to deduce
\begin{equation}\notag
\begin{split}
\int |\nabla^2 \bfu||\nabla^3 \bfu||\nabla^4 \bfu|\, \ud x\aleq& \|\nabla^2\bfu\|_{L^4}\|\nabla^3\bfu\|_{L^4}\|\nabla^4\bfu\|_{L^2}\\
\aleq& \|\nabla\bfu\|_{L^\infty}^{\frac56}\|\nabla^4\bfu\|_{L^2}^{\frac16}\|\nabla\bfu\|_{L^\infty}^{\frac16}\|\nabla^4\bfu\|_{L^2}^{\frac56}\|\nabla^4\bfu\|_{L^2}\\
\aleq& \left(1+\|\bfomg\|_{L^\infty}\log(2+\|\nabla^4\bfu\|_{L^2}) + \nrm{\bfu}_{L^{2}} \right)\|\nabla^4\bfu\|_{L^2}^2.
\end{split}
\end{equation}
Analogous analysis shows that
\begin{equation}\notag
\int \partial^I(\bfu\cdot\nabla\bfB)\cdot\partial^I\bfB\,\ud x=\sum_{I'+I''=3, |I''|\leq 3}c_{I',I''} \int (\partial^{I'}\bfu\cdot\nabla)\partial^{I''}\bfB\cdot\partial^I\bfB \,\ud x
\end{equation}
due to the cancellation
\[ \int (\bfu\cdot\nabla) \partial^I\bfB\cdot\partial^I\bfB\,\ud x =0,\]
and hence
\begin{equation}\notag
\begin{split}
&\left|\int \partial^I(\bfu\cdot\nabla\bfB)\cdot\partial^I\bfB\,\ud x \right|\\
\aleq& \sum_{m_1+m_2=5, 1\leq m_2\leq 4}\int |\nabla^{m_1} \bfu||\nabla^{m_2} \bfB||\nabla^4 \bfB|\, \ud x\\
\aleq& \left(1+\|\bfomg\|_{L^\infty} \log(2+\|\nabla^4\bfu\|_{L^2}) + \|\nb \bfB\|_{L^\infty} + \nrm{\bfu}_{L^{2}} \right)\left(\|\nabla^4\bfu\|_{L^2}^2+\|\nabla^4\bfB\|_{L^2}^2\right) \\
\aleq& \left(1+\|\bfomg\|_{L^\infty} \log(2+\|\nabla^4\bfu\|_{L^2}) +\nrm{\nb(\bfu - \bfJ)}_{L^{\infty} +  \nrm{\bfu}_{L^{2}} + \nrm{\bfB}_{L^{\infty}}}  \right)\left(\|\nabla^4\bfu\|_{L^2}^2+\|\nabla^4\bfB\|_{L^2}^2\right),
\end{split}
\end{equation}
where we used \eqref{eq:cont-crit-H:pf:nbB} in the last inequality. Next, note that
\begin{equation}\notag
\begin{split}
&\int \partial^I(\bfB\cdot\nabla\bfB)\cdot\partial^I\bfu\,\ud x+\int \partial^I(\bfB\cdot\nabla\bfu)\cdot\partial^I\bfB\,\ud x\\
=& \sum_{I'+I''=3, |I''|\leq 3}c_{I',I''} \int (\partial^{I'}\bfB\cdot\nabla)\partial^{I''}\bfB\cdot\partial^I\bfu \,\ud x\\
&+ \sum_{I'+I''=3, |I''|\leq 3}c_{I',I''} \int (\partial^{I'}\bfB\cdot\nabla)\partial^{I''}\bfu\cdot\partial^I\bfB \,\ud x
\end{split}
\end{equation}
thanks to the cancellation
\[\int (\bfB\cdot\nabla) \partial^I\bfB\cdot\partial^I\bfu\,\ud x+\int (\bfB\cdot\nabla) \partial^I\bfu\cdot\partial^I\bfB\,\ud x =0.\]
Again, a similar analysis gives
\begin{equation}\notag
\begin{split}
&\left|\int \partial^I(\bfB\cdot\nabla\bfB)\cdot\partial^I\bfu\,\ud x+\int \partial^I(\bfB\cdot\nabla\bfu)\cdot\partial^I\bfB\,\ud x\right|\\
\aleq& \left(1+\|\bfomg\|_{L^\infty} \log(2+\|\nabla^4\bfu\|_{L^2}) + \|\nb \bfB\|_{L^\infty}+ \nrm{\bfu}_{L^{2}} \right)\left(\|\nabla^4\bfu\|_{L^2}^2+\|\nabla^4\bfB\|_{L^2}^2\right).
\end{split}
\end{equation}
The last integral on the right hand side of \eqref{est-h4-uB-2} is estimated the same way as in the proof for E-MHD, but using \eqref{eq:cont-crit-H:pf:nbB} and \eqref{eq:cont-crit-H:pf:nbnbB} instead of \eqref{eq:cont-crit-J:pf:nbB} and \eqref{eq:cont-crit-J:pf:nbnbB}, respectively. Hence collecting the estimates to be applied to \eqref{est-h4-uB-2} yields
\begin{equation}\notag
\begin{split}
	& \frac{\ud}{\ud t} \left(\nrm{\nb^{4} \bfu}_{L^{2}}^{2}+\nrm{\nb^{4} \bfB}_{L^{2}}^{2} \right) +\nu\nrm{\nb^{5} \bfu}_{L^{2}}^{2} + \nrm{\nb^{5} \bfB}_{L^{2}}^{2} \\
	\aleq& \left(\log (2 + \nrm{\nb^{4} \bfu}_{L^{2}}+ \nrm{\nb^{4} \bfB}_{L^{2}}) + (1+\eta^{-1}) \nrm{\bfB}_{L^{\infty}} + \nrm{\bfu}_{L^{2}} +1 \right)\\
	&\cdot \left(\nrm{\bfomg}_{L^{\infty}}+ \nrm{\nb (u-\bfJ)}_{L^{\infty}} + \nrm{\bfB}_{L^{\infty}}\right)\left(\nrm{\nb^{4} \bfu}_{L^{2}}^{2}+\nrm{\nb^{4} \bfB}_{L^{2}}^{2}\right).
\end{split}
\end{equation}
It is noteworthy that we did not use the term $\nu\nrm{\nb^{5} \bfu}_{L^{2}}^{2}$ in the estimates above and hence the inequality holds for all $\nu\geq0$, including the inviscid case $\nu=0$.
In view of the condition \eqref{eq:cont-crit-H}, it then follows from the generalized Gr\"onwall's inequality that $\nrm{\nb^{4} \bfu(t)}_{L^{2}}^{2}+\nrm{\nb^{4} \bfB(t)}_{L^{2}}^{2}$ is uniformly bounded up to the time $t=T$. \hfill \qedsymbol

\end{document}